\documentclass{amsart}
\usepackage{amscd,amsfonts}
\def\Image{\operatorname{Im}}
\def\Min{\operatorname{Min}}

\def\N{\mathbb{N}}
\def\Z{\mathbb{Z}}
\def\R{\mathbb{R}}

\newtheorem{Theorem}{Theorem}[section]

\newtheorem{Lemma}[Theorem]{Lemma}
\newtheorem{Proposition}[Theorem]{Proposition}
\newtheorem{Thm}{Theorem}

\theoremstyle{definition}
\newtheorem{Definition}[Theorem]{Definition}
\newtheorem{Remark}[Theorem]{Remark}

\begin{document}
\sloppy
\title{
On splitting theorems for CAT(0) spaces and 
compact geodesic spaces of non-positive curvature
}
\author{Tetsuya Hosaka} 
\address{Department of Mathematics, Faculty of Education, 
Utsunomiya University, Utsunomiya, 321-8505, Japan}
\date{May 23, 2010}
\email{hosaka@cc.utsunomiya-u.ac.jp}
%
\begin{abstract}
In this paper, 
we show some splitting theorems for CAT(0) spaces 
on which a product group acts geometrically 
and we obtain a splitting theorem for 
compact geodesic spaces of non-positive curvature.
A CAT(0) group $\Gamma$ is said to be {\it rigid}, 
if $\Gamma$ determines its boundary up to homeomorphisms 
of a CAT(0) space on which $\Gamma$ acts geometrically.
C.~Croke and B.~Kleiner have constructed a non-rigid CAT(0) group.
As an application of the splitting theorems for CAT(0) spaces, 
we obtain that 
if $\Gamma_1$ and $\Gamma_2$ are rigid CAT(0) groups 
then so is $\Gamma_1\times \Gamma_2$.
\end{abstract}
\maketitle

\section{Introduction and preliminaries}

The purpose of this paper is to show some splitting theorems for 
CAT(0) spaces on which a product group acts geometrically 
(i.e.\ properly and cocompactly by isometries).
As applications, we obtain a splitting theorem 
for compact geodesic spaces of non-positive curvature 
and rigidity of the product of rigid CAT(0) groups.

On splitting theorems for spaces of non-positive curvature, 
we can find historical research 
in \cite{BH}, \cite{GW}, \cite{LY}, \cite{S} and \cite{Mo}.

After some preliminaries on CAT(0) spaces and their boundaries in Section~2, 
we first show the following splitting theorem 
by a similar argument to the proof of \cite[Proposition~II.6.23]{BH} in Section~3.

\begin{Thm}\label{Thm1}
Suppose that a group $\Gamma=\Gamma_1\times \Gamma_2$ acts geometrically 
on a CAT(0) space $X$.
If $\Gamma_1$ acts cocompactly on the convex hull $C(\Gamma_1 x_0)$ 
of some $\Gamma_1$-orbit, 
then 
there exists a closed, convex, $\Gamma$-invariant, 
quasi-dense subspace $X'\subset X$ such that 
$X'$ splits as a product $X_1 \times X_2$, 
$\Gamma_1$ naturally acts geometrically on $X_1$, and 
$\Gamma_2$ acts geometrically on $X_2$ by projection.
Here each subspace of the form 
$X_1\times \{x_2\}$ is the closed convex hull of some $\Gamma_1$-orbit.
\end{Thm}

Here we define the induced action by projection in Section~2.

From Theorem~\ref{Thm1}, 
we obtain the following splitting theorem 
which corresponds to \cite[Theorem~II.6.21]{BH}.
(N.~Monod \cite[Corollary~10]{Mo} independently has proved 
the generalization of this splitting theorem.)

\begin{Thm}\label{Thm2}
Suppose that a group $\Gamma=\Gamma_1\times \Gamma_2$ acts geometrically 
on a CAT(0) space $X$.
If the center of $\Gamma$ is finite, then 
there exists a closed, convex, $\Gamma$-invariant, 
quasi-dense subspace $X'\subset X$ such that 
$X'$ splits as a product $X_1 \times X_2$ and 
the action of $\Gamma=\Gamma_1\times \Gamma_2$ on $X'=X_1 \times X_2$ 
is the product action.
\end{Thm}

In \cite{S} and \cite{Mo}, 
some splitting theorems for Hadamard manifolds and CAT(0) spaces have been proved.
In their splitting theorems \cite[Theorem~1]{S} and \cite[Corollary~10]{Mo} and 
also in Theorem~\ref{Thm2}, 
the assumption ``there is not a $\Gamma$-fixed point in $\partial X$'' 
in \cite[Theorem~1]{S} and \cite[Corollary~10]{Mo} 
and the assumption ``the center of $\Gamma$ is finite'' in Theorem~\ref{Thm2} 
are very important, 
and in these cases, the action does split.
Here we note that 
in the case a group $\Gamma$ acts geometrically on a CAT(0) space $X$, 
$\Gamma$ has finite center if and only if 
there is not a $\Gamma$-fixed point in $\partial X$ (cf.\ Lemma~\ref{lem2}).

For example, it is kown that every Coxeter group is a CAT(0) group whose center is finite.

In Section~4, we prove the following splitting theorem 
which is the main theorem in this paper and 
in which we do not need the assumption that ``there is not a $\Gamma$-fixed point in $\partial X$'' 
and that ``the center of $\Gamma$ is finite''.

\begin{Thm}\label{Thm3}
Suppose that a group $\Gamma=\Gamma_1\times \Gamma_2$ acts geometrically on a CAT(0) space $X$.
Then 
there exist finite-index subgroups $G_1$ and $G_2$ of $\Gamma_1$ and $\Gamma_2$ respectively 
and there exist closed convex subspaces $X_1, X_2, X'_1, X'_2$ in $X$ such that 
\begin{enumerate}
\item[(1)] $X_1\times X'_2$ is a closed convex $\Gamma_1\times G_2$-invariant 
quasi-dense subspace of $X$,
\item[(2)] $X'_1\times X_2$ is a closed convex $G_1\times \Gamma_2$-invariant 
quasi-dense subspace of $X$,
\item[(3)] $X'_1$ and $X'_2$ are isometric to 
some quasi-dense subspaces of $X_1$ and $X_2$ respectively,
\item[(4)] there exist geometric actions of $\Gamma_1$ and $\Gamma_2$ 
on $X_1$ and $X_2$ respectively, and
\item[(5)] there exist geometric actions of $G_1$ and $G_2$ 
on $X'_1$ and $X'_2$ respectively.
\end{enumerate}
Hence the boundary $\partial X$ of $X$ is homeomorphic to 
the join $\partial X_1*\partial X_2$ of the boundaries of $X_1$ and $X_2$.
\end{Thm}

In this splitting theorem, 
the center of $\Gamma$ need not be finite 
and also we do not need the assumption 
that there is not a $\Gamma$-fixed point in $\partial X$.
On the other hand, the action does not split in general.
For example, the geometric action of $\Gamma=\Z\times \Z$ on $X=\R\times \R$ 
defined by $(a,b)\cdot (x,y)=(a+x,a+b+y)$ for 
$(a,b)\in\Z\times \Z$ and $(x,y)\in \R\times \R$ does not split 
(of course, the CAT(0) space $X$ splits).

A CAT(0) space $X$ is said to have the {\it geodesic extension property}, 
if every geodesic can be extended to a geodesic line $\R\rightarrow X$.
On CAT(0) spaces with the geodesic extension property, 
we obtain the following splitting theorem as an application of Theorem~\ref{Thm3}.

\begin{Thm}\label{Thm5}
Suppose that a group $\Gamma=\Gamma_1\times \Gamma_2$ acts geometrically 
on a CAT(0) space $X$ with the geodesic extension property.
Then 
$X$ splits as a product $X_1 \times X_2$ and 
there exist geometric actions of $\Gamma_1$ and $\Gamma_2$ 
on $X_1$ and $X_2$, respectively.
Moreover if $\Gamma$ has finite center, 
then $\Gamma$ preserves the splitting, i.e., 
the action of $\Gamma=\Gamma_1\times \Gamma_2$ on $X=X_1 \times X_2$ 
is the product action.
\end{Thm}

In Section~5, 
we study products of CAT(0) groups.
A group $\Gamma$ is called a {\it CAT(0) group}, 
if $\Gamma$ acts geometrically on some CAT(0) space.

Theorem~\ref{Thm3} implies the following.

\begin{Thm}\label{Thm4-0}
$\Gamma_1$ and $\Gamma_2$ are CAT(0) groups if and only if 
$\Gamma_1\times \Gamma_2$ is a CAT(0) group.
\end{Thm}

In \cite{CK}, C.~Croke and B.~Kleiner have proved that 
there exists a CAT(0) group $\Gamma$ and CAT(0) spaces $X$ and $Y$ 
such that $\Gamma$ acts geometrically on $X$ and $Y$ and 
the boundaries of $X$ and $Y$ are not homeomorphic.
On CAT(0) groups and their boundaries, 
it is an important open problem 
when does a CAT(0) group determine its boundary up to homeomorphisms.
A CAT(0) group $\Gamma$ is said to be {\it rigid}, 
if $\Gamma$ determines its boundary up to homeomorphisms 
of a CAT(0) space on which $\Gamma$ acts geometrically.
Then we denote the boundary of the rigid CAT(0) group $\Gamma$ by $\partial \Gamma$.

A conclusion in \cite[Theorem~II.7.1]{BH} and \cite{BR} 
implies that if $\Gamma$ is a rigid CAT(0) group, 
then $\Gamma\times \Z^n$ is also a rigid CAT(0) group for each $n\in\N$.
In \cite{R}, K.~Ruane has proved that 
if $\Gamma_1\times\Gamma_2$ is a CAT(0) group 
and if $\Gamma_1$ and $\Gamma_2$ are hyperbolic groups (in the sense of Gromov) 
then $\Gamma_1\times\Gamma_2$ is rigid.
On rigidity of products of rigid CAT(0) groups, 
we obtain the following theorem from Theorem~\ref{Thm3} 
as a generalization of these results.

\begin{Thm}\label{Thm4}
If $\Gamma_1$ and $\Gamma_2$ are rigid CAT(0) groups, 
then so is $\Gamma_1\times \Gamma_2$, and 
the boundary $\partial(\Gamma_1\times \Gamma_2)$ 
is homeomorphic to the join $\partial\Gamma_1*\partial\Gamma_2$ 
of the boundaries of $\Gamma_1$ and $\Gamma_2$.
\end{Thm}



\section{Lemmas on CAT(0) spaces and their boundaries}

In this section, 
we introduce some properties of CAT(0) spaces and their boundaries 
and show some lemmas needed later.

Definitions and details of CAT(0) spaces and their boundaries 
are found in \cite{BH} and \cite{GH}.
A {\it geometric} action on a CAT(0) space 
is an action by isometries which is proper (\cite[p.131]{BH}) and cocompact.
We note that every CAT(0) space on which some group acts 
geometrically is a proper space (\cite[p.132]{BH}).

Let $X$ be a metric space and let $\gamma$ be an isometry of $X$.
Then the {\it translation length} of $\gamma$ is defined as 
$|\gamma|=\inf\{d(x,\gamma x)\,|\, x\in X\}$, and 
the {\it minimal set} of $\gamma$ is defined as 
$\Min(\gamma)=\{x\in X\,|\, d(x,\gamma x)=|\gamma|\}$.
If $\Gamma$ is a group acting by isometries on $X$, then 
$\Min(\Gamma):=\bigcap_{\gamma\in\Gamma}\Min(\gamma)$.
An isometry $\gamma$ is said to be {\it semi-simple} 
if $\Min(\gamma)$ is non-empty.

Also an isometry $\gamma$ of a metric space $X$ is called 
\begin{enumerate}
\item[(1)] {\it elliptic} if $\gamma$ has a fixed point,
\item[(2)] {\it hyperbolic} if $\gamma$ is semi-simple and not elliptic, and 
\item[(3)] {\it parabolic} if $\gamma$ is not semi-simple.
\end{enumerate}

Now we define induced actions by projection.

\begin{Definition}
Suppose that a group $\Gamma=\Gamma_1\times \Gamma_2$ acts geometrically 
on a CAT(0) space $X=X_1\times X_2$.
In general, the action does not split.
Let $p_1:X_1\times X_2\rightarrow X_1$ be the projection.
We can define an action ``$*$'' of $\Gamma_1$ on $X_1$ 
by $\gamma_1*x_1=p_1(\gamma_1\cdot x_1)$ 
for $\gamma_1\in \Gamma_1$ and $x_1\in X_1$ (cf.\ \cite{BR}), 
where ``$\cdot$'' is the original action of $\Gamma$ on $X$.
Then we say that 
the action ``$*$'' of $\Gamma_1$ on $X_1$ is the induced action by the projection of 
the action ``$\cdot$''.

In this paper, we often say that 
{\it $\Gamma_1$ acts geometrically on $X_1$ by projection}, 
if the induced action of $\Gamma_1$ on $X_1$ by the projection of 
the original action of $\Gamma$ on $X$ is geometric.

Also, in this paper, we often say that 
{\it $\Gamma_1$ acts geometrically on $X_1$ by restriction}, 
if $\Gamma_1$ acts geometrically on $X_1$ 
by the restriction of the original action of $\Gamma$ on $X$ 
(hence, in such a case, $X_1$ is $\Gamma_1$-invariant).
\end{Definition}

The following theorem is called the Flat Torus Theorem 
(cf.\ \cite{AB}, \cite[Proposition~1.1]{BR}, \cite[Theorem~II.7.1]{BH}).

\begin{Theorem}[{\cite[Theorem~II.7.1]{BH}}]\label{thmBR}
Let $G$ be a group and let $A$ be a free abelian group of rank $n$.
Suppose that $\Gamma=G\times A$ acts geometrically on a CAT(0) space $X$.
Then 
\begin{enumerate}
\item[(1)] $\Min(A)=\bigcap_{\alpha\in A}\Min(\alpha)$ is 
a closed, convex, $\Gamma$-invariant, quasi-dense subspace of $X$ that 
splits as a product $Y\times Z$, where $Z$ is isometric to $\R^n$,
\item[(2)] $G$ acts geometrically on $Y$ by projection, and 
\item[(3)] $A$ acts geometrically on $Z$ by restriction 
(moreover, $Z$ is the convex hull of some orbit $Ax_0$ of $X$).
\end{enumerate}
\end{Theorem}

Here a subset $X'$ of a metric space $X$ is said to be {\it quasi-dense} if 
there exists a number $N>0$ such that each point of $X$ is $N$-close 
to some point of $X'$, i.e., $B(X',N)=X$.

Some results on ``$\Gamma$-invariant quasi-dense subspace'' in this paper 
relate to \cite[Corollary~2.7]{AB} and \cite{BS}.

We show some lemmas needed later.

\begin{Lemma}\label{lem1}
Let $\gamma$ be a parabolic isometry of 
a proper cocompact CAT(0) space $X$.
Then there exists a non-fixed point $\alpha\in\partial X$ of $\gamma$.
\end{Lemma}

\begin{proof}
Let $\gamma$ be a parabolic isometry of a proper cocompact CAT(0) space $X$.
Then the isometry $\gamma$ of $X$ induces a homeomorphism of 
the boundary $\partial X$ (cf. \cite[Corollary~II.8.9]{BH}).
We suppose that 
the induced homeomorphism 
$\gamma:\partial X \rightarrow \partial X$ is the identity.

Since $\gamma$ is a parabolic isometry, 
$\Min(\gamma)$ is the empty set.
Hence $X$ is non-compact and $\partial X$ is non-empty.
This means that there exists a geodesic ray in $X$.
Since $X$ is cocompact, 
there exists a geodesic line $\sigma_1:\R\rightarrow X$ 
by \cite[Lemma~II.9.34]{BH}.
Let $Y_1$ be the union of the images of all geodesic lines parallel to $\sigma_1$.
For each geodesic line $\tau$ parallel to $\sigma_1$, 
by hypothesis, $\gamma\tau(\infty)=\tau(\infty)$ and 
$\gamma\tau(-\infty)=\tau(-\infty)$.
Hence $\gamma\tau$ is also parallel to $\sigma_1$ 
by the Flat Strip Theorem \cite[Theorem~II.2.13]{BH}. 
Thus $Y_1$ is $\gamma$-invariant.
By \cite[Theorem~II.2.14]{BH}, 
$Y_1$ is isometric to $X_1 \times \R$ for some convex subspace $X_1\subset Y_1$.
Then the restriction of $\gamma$ to $Y_1=X_1 \times \R$ 
splits as $(\gamma_1,\gamma_1')$, 
where $\gamma_1$ and $\gamma_1'$ are isometries of $X_1$ and $\R$, respectively.
Since $\gamma'_1$ is semi-simple (cf.\ \cite[Proposition~II.6.5]{BH}) 
and $\gamma$ is not semi-simple, 
by \cite[Proposition~II.6.9]{BH}, 
$\gamma_1$ is not semi-simple, i.e., 
$\gamma_1$ is a parabolic isometry.
We show that 
the induced homeomorphism 
$\gamma_1:\partial X_1\rightarrow \partial X_1$ is the identity.
Let $\xi$ be a geodesic ray in $X_1$.
Then $\gamma \xi$ is a geodesic ray in $\gamma X_1$ and 
$\gamma \xi(\infty)=\xi(\infty)$ 
because $\gamma$ is the identity of $\partial X$.
Hence $\xi$ and $\gamma\xi$ are parallel and 
$\gamma_1\xi(\infty)=\xi(\infty)$ 
by the definition of $\gamma_1$ 
(since $\gamma_1$ is the projection on $X_1$ of $\gamma$).
Thus $\gamma_1:\partial X_1\rightarrow \partial X_1$ is the identity.

Since $\gamma_1$ is a parabolic isometry, 
$X_1$ is non-compact and $\partial X_1$ is non-empty, hence 
there exists a geodesic ray $\xi_1$ in $X_1$.
Then 
$$ X\supset Y_1=X_1\times \R \supset \Image \xi_1 \times \R. $$
Since $X$ is cocompact, 
there exists a subspace $Z$ of $X$ which is isometric to $\R^2$ 
by \cite[Lemma~II.9.34]{BH}.
Here we redefine $\sigma_1$ as a geodesic line in $Z$, 
and we also redefine $Y_1$, $X_1$ and $\gamma_1$ from $\sigma_1$.
Then there exists a geodesic line $\sigma_2$ in $X_1$.
Indeed we can construct $\sigma_2$ as $\Image \sigma_2=Z\cap X_1$.
Let $Y_2$ be the union of the images of 
all geodesic lines parallel to $\sigma_2$ in $X_1$.
Here $Y_2$ is $\gamma_1$-invariant and 
isometric to $X_2 \times \R$ for some convex subspace $X_2$ of $Y_2$ 
by \cite[Theorem~II.2.14]{BH}.
Then the restriction of $\gamma_1$ to $Y_2=X_2 \times \R$ 
splits as $(\gamma_2,\gamma'_2)$, 
$\gamma_2$ is a parabolic isometry of $X_2$ and 
the induced homeomorphism $\gamma_2$ of $\partial X_2$ is the identity.

By iterating this argument, 
for any $n\in \N$, 
there exists a sequence 
\begin{align*} 
X &\supset Y_1=X_1\times \R \supset Y_2\times \R=X_2\times \R^2 \\
&\supset Y_3\times \R^2=X_3\times \R^3\supset \cdots \supset X_n\times \R^n,
\end{align*}
that is, $\R^n\subset X$ for any $n\in \N$.
However this contradicts \cite[Lemma~II.7.4]{BH}.
Thus $\gamma:\partial X \rightarrow \partial X$ is not the identity, 
i.e., $\gamma \alpha \neq \alpha$ for some $\alpha\in\partial X$.
\end{proof}

\begin{Lemma}\label{lem2}
Suppose that a group $\Gamma$ acts geometrically on a CAT(0) space $X$.
The center of $\Gamma$ is finite if and only if 
there does not exist $\alpha\in\partial X$ such that $\Gamma \alpha=\alpha$.
\end{Lemma}

\begin{proof}
We first show that 
if the center of $\Gamma$ is finite then 
there does not exist $\alpha\in\partial X$ such that $\Gamma \alpha=\alpha$.

Suppose that there exists 
$\alpha\in\partial X$ such that $\Gamma \alpha=\alpha$.
Then we show that $\Gamma$ has infinite center.
Let $x_0\in X$ 
and let $\xi:[0,\infty)\rightarrow X$ be the geodesic ray 
such that $\xi(0)=x_0$ and $\xi(\infty)=\alpha$.
We note that 
$\delta\xi(\infty)=\xi(\infty)$ for any $\delta\in\Gamma$ 
because $\Gamma \alpha=\alpha$.
Since the action of $\Gamma$ on $X$ is cocompact, 
$\Gamma B(x_0,N)=X$ for some $N>0$.
Then there exists a sequence $\{\gamma_i\,|\,i\in\N\}\subset \Gamma$ 
such that $d(\gamma_i x_0,\Image \xi)\le N$ 
and $\gamma_i\neq \gamma_j$ if $i\neq j$.

Now we show that 
for each $\delta\in\Gamma$, 
there exists a subsequence 
$\{\gamma_{i_j}\,|\,j\in\N\}\subset\{\gamma_i\,|\,i\in\N\}$ 
such that 
$(\gamma_{i_1}\gamma_{i_j}^{-1})\delta=\delta(\gamma_{i_1}\gamma_{i_j}^{-1})$
for any $j\in \N$.

Let $\delta\in\Gamma$.
Since $\delta\xi(\infty)=\xi(\infty)$, 
there exists a number $K_\delta>0$ such that 
$d(\xi(r),\delta\xi(r))\le K_\delta$ for any $r\ge 0$.
For each $i\in\N$, 
$d(\gamma_i x_0,\xi(r_i))\le N$ for some $r_i\ge 0$ 
because $d(\gamma_i x_0,\Image \xi)\le N$.
Then 
\begin{align*}
d(x_0,\gamma_i^{-1}\delta\gamma_i x_0)&=d(\gamma_ix_0,\delta\gamma_i x_0) \\
&\le d(\gamma_i x_0,\xi(r_i))+d(\xi(r_i),\delta\xi(r_i))
+d(\delta\xi(r_i),\delta\gamma_i x_0) \\
&\le N+K_\delta +d(\xi(r_i),\gamma_i x_0) \\
&\le 2N+K_\delta.
\end{align*}
Since the action of $\Gamma$ on $X$ is proper, 
$\{\gamma\in\Gamma\,|\, d(x_0,\gamma x_0)\le 2N+K_\delta\}$ 
is a finite set.
Hence 
there exists a subsequence 
$\{\gamma_{i_j}\,|\,j\in\N\}\subset\{\gamma_i\,|\,i\in\N\}$ 
such that $\gamma_{i_j}^{-1}\delta\gamma_{i_j}=\gamma_{i_1}^{-1}\delta\gamma_{i_1}$ 
for any $j\in \N$.
Then 
$(\gamma_{i_1}\gamma_{i_j}^{-1})\delta=\delta(\gamma_{i_1}\gamma_{i_j}^{-1})$
for any $j\in \N$.

The CAT(0) group $\Gamma$ is finitely presented (cf.\ \cite[Theorem~III.1.1(1)]{BH}).
Let $\{\delta_1,\dots,\delta_n\}$ be a generating set of $\Gamma$.
Then by iterating the above argument, 
there exists a subsequence 
$\{g_i\,|\,i\in\N\}\subset\{\gamma_i\,|\,i\in\N\}$ 
such that $(g_1g_i^{-1})\delta_k=\delta_k(g_1g_i^{-1})$ 
for any $i\in\N$ and $k\in\{1,\dots,n\}$.
Since $\{\delta_1,\dots,\delta_n\}$ generates $\Gamma$, 
$$(g_1g_i^{-1})\delta=\delta(g_1g_i^{-1})$$ 
for any $i\in\N$ and $\delta\in \Gamma$.
Here we note that $g_1g_i^{-1}\neq g_1g_j^{-1}$ if $i\neq j$ 
by the construction of $\{\gamma_i\,|\,i\in\N\}$.
Thus the center of $\Gamma$ contains 
the infinite set $\{g_1g_i^{-1}\,|\,i\in\N\}$.

Next we show that 
if there does not exist $\alpha\in\partial X$ such that $\Gamma \alpha=\alpha$, 
then the center of $\Gamma$ is finite.

Suppose that the center $Z(\Gamma)$ of $\Gamma$ is infinite.
Since the set $\{\delta x_0\,|\,\delta\in Z(\Gamma)\}$ is unbounded in $X$, 
there exists a sequence 
$$ \{\delta_i x_0 \,|\,i\in\N \} \subset \{\delta x_0\,|\,\delta\in Z(\Gamma)\}$$ 
which converges to some point $\alpha_0\in \partial X$ in $X\cup\partial X$.
Let $\gamma\in\Gamma$.
Then the sequence $\{\gamma\delta_i x_0\,|\,i\in\N \}$ converges to $\gamma\alpha_0$.
Here since $\gamma\delta_i=\delta_i\gamma$,
$$ \{\gamma\delta_ix_0\,|\,i\in\N \}=\{\delta_i\gamma x_0\,|\,i\in\N \}=\{\delta_i y_0\,|\,i\in\N \},$$
where $y_0=\gamma x_0$, 
which converges to $\alpha_0$ in $X\cup\partial X$.
Hence $\gamma\alpha_0=\alpha_0$ for any $\gamma\in\Gamma$.
Thus $\Gamma \alpha_0=\alpha_0$.
\end{proof}

\section{Some splitting theorems}

We first generalize the splitting theorem \cite[Proposition~II.6.23]{BH} as follows.
This plays a key role to prove the main results.

\begin{Theorem}\label{thm1}
Suppose that a group $\Gamma=\Gamma_1\times \Gamma_2$ acts geometrically 
on a CAT(0) space $X$.
If $\Gamma_1$ acts cocompactly on the convex hull $C(\Gamma_1 x_0)$ 
of some $\Gamma_1$-orbit, 
then 
there exists a closed, convex, $\Gamma$-invariant, 
quasi-dense subspace $X'\subset X$ such that 
\begin{enumerate}
\item[(1)] $X'$ splits as a product $X_1 \times X_2$, 
\item[(2)] $\Gamma_1$ acts geometrically on $X_1$ by restriction 
(moreover, each subspace of the form 
$X_1\times \{x_2\}$ is the closed convex hull of some $\Gamma_1$-orbit), and 
\item[(3)] $\Gamma_2$ acts geometrically on $X_2$ by projection.
\end{enumerate}
\end{Theorem}

\begin{proof}
We show this theorem by a similar proof to the one in \cite[pp.240--241]{BH}.

Let $\Sigma$ be the set of closed, convex, non-empty, 
$\Gamma_1$-invariant subspaces of $X$, and 
let $\mathcal{N}$ be the subset of $\Sigma$ consisting of those subspaces 
which are minimal with respect to inclusion.
We note that the member of $\mathcal{N}$ are disjoint and 
each member of $\mathcal{N}$ is the closed convex hull $C(\Gamma_1 x)$ of 
some $\Gamma_1$-orbit.

We mainly check corresponding claims in \cite[pp.240--241]{BH}.
From Zorn's lemma, we first obtain Claim~1: $\mathcal{N}$ is non-empty.

For $C_1,C_2\in \mathcal{N}$, 
let $p_i:X\rightarrow C_i$ denote the projection of $X$ onto $C_i$ 
and let $d=d(C_1,C_2):=\inf\{d(x_1,x_2)\,|\,x_1\in C_1, x_2\in C_2\}$.
By \cite[p.240]{BH}, 
we also obtain Claim~2:
There exists a unique isometry $j$ of $C_1\times [0,d]$ 
onto the convex hull of $C_1\cup C_2$ such that 
$j(x,0)=x$ and $j(x,d)=p_2(x)$ for each $x\in C_1$.

Let $X'=\bigcup \mathcal{N}$.
Then we check 
Claim~3: $X'$ is a closed, convex, $\Gamma$-invariant, quasi-dense subspace of $X$.
We note that 
for each $x\in X'$ there exists a unique member $C\in \mathcal{N}$ 
such that $x\in C$ 
(indeed $C$ is the convex hull of the $\Gamma_1$-orbit $\Gamma_1 x$).
By \cite[p.240]{BH}, we see that 
$X'$ is a closed, convex, $\Gamma$-invariant subspace of $X$.
We show that $X'$ is quasi-dense in $X$.
Since $\Gamma$ acts cocompactly on $X$, 
$\Gamma B(x_0,N) =X$ for some $x_0\in X'$ and $N>0$.
Here $\Gamma x_0\subset X'$ because $X'$ is $\Gamma$-invariant.
Hence 
$$ B(X',N)\supset B(\Gamma x_0,N)=\Gamma B(x_0,N)=X.$$
Thus $X'$ is a quasi-dense subspace of $X$.

Next we check Claim~4:
If $C_1,C_2,C_3\in \mathcal{N}$ and 
if $p_i:X\rightarrow C_i$ denotes the projection onto $C_i$, 
then $p_1=p_1p_3$ on $C_2$.

Let $f=p_1p_3p_2:C_1\rightarrow C_2 \rightarrow C_3 \rightarrow C_1$.
Then $f$ is an isometry of $C_1$ by Claim~2.
To show that $p_1=p_1p_3$ on $C_2$, 
we prove that $f$ is the identity.

We first show that $f$ is an elliptic isometry.
Let $K=d(C_1,C_2)+d(C_2,C_3)+d(C_3,C_1)$.
Then $d(x,f(x))\le K$ for any $x\in C_1$.
Hence for each geodesic ray $\xi$ in $C_1$, 
$d(\xi(r),f\xi(r))\le K$ for any $r\in [0,\infty)$.
This means that 
every geodesic ray $\xi$ in $C_1$ is asymptotic to 
the geodesic ray $f\xi$.
Hence the induced homeomorphism 
$f:\partial C_1\rightarrow \partial C_1$ is the identity.
Thus $f$ is not parabolic isometry by Lemma~\ref{lem1}.
Now we suppose that $f$ is a hyperbolic isometry of $C_1$.
Then there exists an axis $\sigma_1$ of $f$ in $C_1$.
Let $\sigma_2=p_2\sigma_1$ and let $\sigma_3=p_3\sigma_2$.
Here we note that 
$\sigma_2$ and $\sigma_3$ are geodesic lines in $C_2$ and $C_3$ respectively, 
and the geodesic lines 
$\sigma_1$, $\sigma_2$ and $\sigma_3$ are parallel by Claim~2.
Then 
$$p_1\sigma_3=p_1p_3\sigma_2=p_1p_3p_2\sigma_1=f\sigma_1,$$
and $f(\Image\sigma_1)=\Image\sigma_1$.
By \cite[Lemma~II.2.15]{BH}, 
$f$ is the identity on $\Image \sigma_1$. 
This is a contradiction, since $\sigma_1$ is an axis of $f$.
Hence $f$ is not a hyperbolic isometry.
Thus $f$ is elliptic.

We show that $f$ is the identity.
Let $C'_1=\Min (f)$.
Here $C'_1$ is the fixed point set of $f$.
To show that $f$ is the identity, 
we prove that $C'_1=C_1$.
If $C'_1$ is a closed, convex, non-empty, 
$\Gamma_1$-invariant subspace of $X$, 
then $C'_1=C_1$ by the minimality of $C_1$.
Now $C'_1$ is a closed, convex, non-empty subset.
Hence it is sufficient to show that $C'_1$ is $\Gamma_1$-invariant.
Let $\gamma_1\in\Gamma_1$ and let $x\in C'_1$.
Then $f(x)=x$.
We show that $\gamma_1 x\in C'_1$ (i.e.\ $f(\gamma_1 x)=\gamma_1 x$).
Since $p_i$ is $\Gamma_1$-invariant for each $i=1,2,3$, 
\begin{align*}
f(\gamma_1 x)&=p_1p_3p_2(\gamma_1 x)=p_1p_3(\gamma_1 p_2(x))=p_1(\gamma_1 p_3p_2(x)) \\
&=\gamma_1(p_1p_3p_2(x))=\gamma_1 f(x)=\gamma_1 x.
\end{align*}
This means that $\gamma_1 x\in C'_1$.
Hence $C'_1$ is $\Gamma_1$-invariant and 
$C'_1=C_1$, i.e., $f$ is the identity.

\vspace{1mm}

We fix $X_1\in \mathcal{N}$ and 
let $p:X'\rightarrow X_1$ be the orthogonal projection.
Let $X_2$ denote the metric space $(\mathcal{N},d)$, 
where $d(C,C')=\inf\{d(x,x')\,|\,x\in C, x'\in C'\}$.
For each $x\in X'$, 
there exists a unique member $C_x\in \mathcal{N}$ such that $x\in C_x$.
Then by \cite[p.241]{BH}, we obtain 
Claim~5: The map $X'\rightarrow X_1\times X_2$ given by 
$x \mapsto (p(x),C_x)$ is a $\Gamma$-equivariant isometry.

Here $X_1$ is $\Gamma_1$-invariant and 
$\Gamma_1$ acts geometrically on $X_1$.
Finally we show that $\Gamma_2$ acts geometrically on $X_2$ by $\gamma_2 C_x=C_{\gamma_2 x}$ 
(i.e., $\Gamma_2$ acts geometrically on $X_2$ by projection).
It is obvious that the action of $\Gamma_2$ on $X_2$ is cocompact and isometry.
Suppose that the action of $\Gamma_2$ on $X_2$ is not proper.
Then there exists a number $K>0$ such that 
\begin{align*}
&\{\gamma\in\Gamma_2\,|\,B(C_{x_0},K)\cap \gamma B(C_{x_0},K)\neq\emptyset\} \\
&\ =\{\gamma\in\Gamma_2\,|\,B(C_{x_0},K)\cap B(C_{\gamma x_0},K)\neq\emptyset\} 
\end{align*}
is infinite.
Let $\{\gamma_i\,|\,i\in\N\}=
\{\gamma\in\Gamma_2\,|\,B(C_{x_0},K)\cap B(C_{\gamma x_0},K)\neq\emptyset\}$, 
where $\gamma_i\neq \gamma_j$ if $i\neq j$.
Since $d(C_{x_0},C_{\gamma_i x_0})\le 2K$, 
there exists a point $y_0\in C_{x_0}$ such that 
$d(y_0,\gamma_i x_0)\le 2K$ by Claim~2.
Now $\Gamma_1$ acts cocompactly on $C_{x_0}$ 
(which is the closed convex hull of the orbit $\Gamma_1 x_0$), 
and $C_{x_0}\subset\Gamma_1 B(x_0,N)$ for some $N>0$.
There exists $\delta_i\in\Gamma_1$ such that 
$d(\delta_i x_0,y_0)\le N$. 
Then 
$$d(\delta_i x_0,\gamma_i x_0)
\le d(\delta_i x_0,y_0)+d(y_0,\gamma_i x_0)\le N+2K. $$
Thus $d(x_0,\delta_i^{-1}\gamma_i x_0)\le N+2K$ for any $i\in\N$.
Here we note that 
$\delta_i^{-1}\gamma_i\neq \delta_j^{-1}\gamma_j$ if $i\neq j$, 
because $\delta_i,\delta_j\in\Gamma_1$, 
$\gamma_i,\gamma_j\in\Gamma_2$ 
and $\gamma_i\neq \gamma_j$.
This means that the action of $\Gamma=\Gamma_1\times\Gamma_2$ on $X$ 
is not proper, which is a contradiction.
Thus the action of $\Gamma_2$ on $X_2$ is proper, hence, geometric.
\end{proof}

\begin{Remark}\label{Rem3.2}
We suppose that a group $\Gamma=\Gamma_1\times \Gamma_2$ acts geometrically 
on a CAT(0) space $X$ and we suppose that {\it $\Gamma_2$ has finite center}.
Then \cite[Lemma~II.6.24]{BH} implies that 
$\Gamma_1$ acts cocompactly on the convex hull $C(\Gamma_1 x_0)$ 
of some $\Gamma_1$-orbit, that is, 
the condition in Theorem~\ref{thm1} holds.
In this paper, 
we mainly use Theorem~\ref{thm1} under the condition 
that $\Gamma_2$ has finite center.
\end{Remark}

From Theorem~\ref{thm1} and Lemma~\ref{lem2}, 
we obtain the following splitting theorem.

\begin{Theorem}\label{thm2}
Suppose that a group $\Gamma=\Gamma_1\times \Gamma_2$ acts geometrically 
on a CAT(0) space $X$.
If $\Gamma$ has finite center, then 
there exists a closed, convex, 
$\Gamma$-invariant, quasi-dense subspace $X'\subset X$ such that 
$X'$ splits as a product $X_1 \times X_2$ and 
the action of $\Gamma=\Gamma_1\times \Gamma_2$ on $X'=X_1 \times X_2$ 
is the product action.
\end{Theorem}

\begin{proof}
Suppose that $\Gamma$ has finite center.
Then the centers of $\Gamma_1$ and $\Gamma_2$ are finite.
Let $X'$ be a minimal closed convex non-empty $\Gamma$-invariant subspace of $X$
(there exists such $X'$ by Zorn's lemma).
Then for each $x\in X'$, 
the convex hull $C(\Gamma x)$ of the $\Gamma$-orbit $\Gamma x$ 
is $X'$ by the minimality of $X'$.
Let $x_0\in X'$ and let $X_1=C(\Gamma_1 x_0)$.
Since the action of $\Gamma_1$ on $X_1$ is cocompact by \cite[Lemma~II.6.24]{BH}, 
$X'$ splits as $X_1\times X_2$ by Theorem~\ref{thm1}.
Then we show that $X_2=C(\Gamma_2 x_0)$ 
and the action of $\Gamma_2$ on $X_1$ is trivial.

Let $p:X'\rightarrow X_1$ be the projection.
For each $\gamma\in\Gamma_2$ and $x\in X_1$, 
we define $\gamma*x=p(\gamma x)$.
Then ``$*$'' is an isometry action of $\Gamma_2$ on $X_1$.
Here
for each $\gamma_1\in \Gamma_1$, $\gamma_2\in \Gamma_2$ and $x\in X_1$, 
\begin{align*}
\gamma_1(\gamma_2*x)&=\gamma_1p(\gamma_2 x)=p(\gamma_1\gamma_2 x) \\
&=p(\gamma_2\gamma_1 x)=\gamma_2*(\gamma_1 x).
\end{align*}
Hence $\gamma_1(\gamma_2*x)=\gamma_2*(\gamma_1 x)$.

Now we show that the orbit $\Gamma_2 *x_0$ is a bounded set.
Suppose that $\Gamma_2 *x_0$ is unbounded. 
Then there exists a sequence $\{g_i*x_0\}\subset \Gamma_2 *x_0$ 
which converges to some point $\alpha\in\partial X_1$ in $X_1\cup\partial X_1$.
For each $\delta\in\Gamma_1$, 
$\{\delta(g_i*x_0)\}=\{g_i*(\delta x_0)\}$ converges to $\alpha$.
On the other hand, 
$\{\delta(g_i*x_0)\}$ converges to $\delta\alpha$.
Hence $\delta\alpha=\alpha$ for any $\delta\in\Gamma_1$.
This contradicts Lemma~\ref{lem2}, 
since the center of $\Gamma_1$ is finite.
Thus $\Gamma_2* x_0$ is a bounded set in $X_1$.
By \cite[Corollary~II.2.8(1)]{BH}, 
there exists a point $y_0\in X_1$ such that $\Gamma_2*y_0=y_0$.
Here we may retake the basepoint $x_0$ as $y_0$.
Then $X_1=C(\Gamma_1 x_0)$ and $X_2=C(\Gamma_2 x_0)$.
By the construction of $X_1\times X_2$ in the proof of Theorem~\ref{thm1}, 
the action of $\Gamma_2$ on $X_1$ is trivial.
\end{proof}

A metric space $Y$ is said to be of {\it non-positive curvature}, 
if $Y$ is a locally CAT(0) space (cf.\ \cite[p.159]{BH}).
Let $Y$ be a compact geodesic space of non-positive curvature.
Then the universal covering $X$ of $Y$ is a CAT(0) space 
by the Cartan-Hadamard theorem, 
and we can think of $Y$ as the quotient $\Gamma\backslash X$ of $X$,
where $\Gamma$ is the fundamental group of $Y$ acting freely 
and properly by isometries on $X$ (cf.\ \cite[p.237]{BH}).

As an application of Theorem~\ref{thm2}, 
we can obtain the following splitting theorem for 
compact geodesic spaces of non-positive curvature (i.e.\ locally CAT(0) spaces).
This splitting theorem is relate to the splitting theorems in 
\cite[Theorem~II.6.22]{BH}, \cite{GW}, \cite{LY} and \cite{S}.

\begin{Theorem}\label{thm6}
Let $Y$ be a compact geodesic space of non-positive curvature (i.e.\ locally CAT(0)).
Suppose that the fundamental group of $Y$ splits as a product 
$\Gamma=\Gamma_1\times \Gamma_2$ 
and that $\Gamma$ has trivial center.
Then there exists a deformation retract $Y'$ of $Y$ 
which splits as a product $Y_1\times Y_2$ 
such that the fundamental group of $Y_i$ is $\Gamma_i$ for each $i=1,2$.
\end{Theorem}

\begin{proof}
Let $X$ be the universal covering of $Y$.
Then $X$ is a CAT(0) space, $\Gamma$ acts freely and properly by isometries on $X$, 
and $\Gamma\backslash X=Y$.
Since $\Gamma=\Gamma_1\times \Gamma_2$ acts geometrically on $X$ and 
$\Gamma$ has finite center, 
by Theorem~\ref{thm2}, 
there exists a closed, convex, 
$\Gamma$-invariant, quasi-dense subspace $X'\subset X$ such that 
$X'$ splits as a product $X_1 \times X_2$ and 
the action of $\Gamma=\Gamma_1\times \Gamma_2$ on $X'=X_1 \times X_2$ 
is the product action.
Let $Y_i=\Gamma_i\backslash X_i$ for each $i=1,2$.
Then
\begin{align*} 
Y&=\Gamma\backslash X\supset \Gamma\backslash X' 
=(\Gamma_1\times\Gamma_2)\backslash(X_1\times X_2) \\
&=\Gamma_1\backslash X_1\times \Gamma_2\backslash X_2=Y_1\times Y_2, 
\end{align*}
and $\pi_1Y_i=\Gamma_i$ for each $i=1,2$, 
since $\Gamma_i$ acts freely and properly by isometries on $X_i$.

We show that $Y_1\times Y_2$ is a deformation retract of $Y$.
Since $X'=X_1\times X_2$ is a quasi-dense subspace of $X$, 
there exists $N\ge 0$ such that $B(X',N)=X$.
For each $t\in[0,N]$, 
let $p_t:X\rightarrow B(X',t)$ be the orthogonal projection.
Then the map $F:X\times [0,N]\rightarrow X$ 
defined by $F(x,t)=p_{N-t}(x)$ is a deformation retraction of $X'$ onto $X$.
Since $X'$ is $\Gamma$-invariant, 
$B(X',t)$ is also $\Gamma$-invariant 
and the square
$$ \begin{CD}
X @>{p_t}>> B(X',t) \\
@V{\gamma}VV @V{\gamma}VV \\
X @>{p_t}>> B(X',t) 
\end{CD} $$
is commutative for each $t\in[0,N]$ and $\gamma\in\Gamma$.
Hence $F$ induces a deformation retract 
$$F':\Gamma\backslash X \times [0,N]\rightarrow \Gamma\backslash X$$
of $\Gamma\backslash X'=Y_1\times Y_2$ onto $\Gamma\backslash X=Y$.
\end{proof}

\section{A splitting theorem for CAT(0) spaces 
on which a product group acts geometrically}

In this section, 
we show a splitting theorem in the case a CAT(0) group has infinite center.

We first define an essential subgroup of a group.

\begin{Definition}
Let $\Gamma$ be a group and let $G$ be a subgroup of $\Gamma$.
In this paper, we say that $G$ is an {\it essential subgroup} of $\Gamma$, 
if $G$ is a finite-index subgroup of $\Gamma$ 
which splits as a product $\Gamma'\times A$, where 
$\Gamma'$ has finite center and 
$A$ is isomorphic to $\Z^n$ for some $n\ge 0$.
\end{Definition}

The following lemma is known.

\begin{Lemma}[{\cite[Lemma~2.1]{Ho}}]\label{lem5.1}
If a group $\Gamma_1\times \Gamma_2$ acts geometrically on some CAT(0) space $X$, 
then there exists an essential subgroup of $\Gamma_i$ for each $i=1,2$.
\end{Lemma}

The purpose of this section is to prove the following main splitting theorem.

\begin{Theorem}\label{Thm}
Suppose that a group $\Gamma=\Gamma_1\times \Gamma_2$ acts geometrically on a CAT(0) space $X$.
Then for essential subgroups $G_1$ and $G_2$ of $\Gamma_1$ and $\Gamma_2$ respectively, 
there exist closed convex subspaces $X_1, X_2, X'_1, X'_2$ in $X$ such that 
\begin{enumerate}
\item[(1)] $X_1\times X'_2$ is a closed convex $\Gamma_1\times G_2$-invariant 
quasi-dense subspace of $X$,
\item[(2)] $X'_1\times X_2$ is a closed convex $G_1\times \Gamma_2$-invariant 
quasi-dense subspace of $X$,
\item[(3)] $X'_1$ and $X'_2$ are isometric to 
some quasi-dense subspaces of $X_1$ and $X_2$ respectively,
\item[(4)] there exist geometric actions of $\Gamma_1$ and $\Gamma_2$ 
on $X_1$ and $X_2$ respectively, and
\item[(5)] there exist geometric actions of $G_1$ and $G_2$ 
on $X'_1$ and $X'_2$ respectively.
\end{enumerate}
Hence the boundary $\partial X$ of $X$ is homeomorphic to 
the join $\partial X_1*\partial X_2$ of the boundaries of $X_1$ and $X_2$.
\end{Theorem}

We can obtain the following proposition from 
Theorems~\ref{thmBR} and \ref{thm1}, Remark~\ref{Rem3.2} and Lemma~\ref{lem5.1} 
by a similar argument to the proof of the splitting theorem in \cite{Ho}.
To prove Theorem~\ref{Thm}, we investigate details of the argument.

\begin{Proposition}\label{prop6.1}
Suppose that a group $\Gamma=\Gamma_1\times \Gamma_2$ acts geometrically 
on a CAT(0) space $X$.
Then for an essential subgroup $G_2$ of $\Gamma_2$, 
there exists a closed, convex, $\Gamma_1\times G_2$-invariant, 
quasi-dense subspace $X'$ of $X$ such that 
$X'$ splits as a product $X_1 \times X'_2$ and 
there exist geometric actions of $\Gamma_1$ and $G_2$ on $X_1$ and $X'_2$ respectively.
\end{Proposition}

We provide two remarks needed in the proof of Theorem~\ref{Thm}.

\begin{Remark}\label{rem1}
Suppose that 
a group $\Gamma=\Gamma_1\times\Gamma_2\times A$ 
acts geometrically on a CAT(0) space $X$, 
where the centers of $\Gamma_1$ and $\Gamma_2$ are finite and 
$A$ is isomorphic to $\Z^n$ for some $n$.

By Theorem~\ref{thmBR}, 
$\Min(A)$ splits as a product $Z\times Y$ 
which is a quasi-dense subspace of $X$, 
where $Z$ is isometric to $\R^n$, 
$A$ acts geometrically on $Z$ by restriction and 
$\Gamma_1\times\Gamma_2$ acts geometrically on $Y$ by projection.
Since $\Gamma_1\times\Gamma_2$ has finite center, by Theorem~\ref{thm2}, 
$Y$ contains a quasi-dense subspace $Y_1\times Y_2$, 
where $\Gamma_1$ and $\Gamma_2$ 
acts geometrically on $Y_1$ and $Y_2$ respectively 
by the restriction of the action 
of $\Gamma_1\times\Gamma_2$ on $Y$ by projection.
Here 
$$ X\supset Z\times Y \supset Z\times Y_1\times Y_2.$$
We note that for the original action of $\Gamma$ on $X$, 
$A$ acts geometrically on $Z$ by restriction and 
$\Gamma_i$ acts geometrically on $Y_i$ by projection for each $i=1,2$.

On the other hand, 
since $\Gamma_1\times(\Gamma_2\times A)$ acts geometrically on $X$ 
and $\Gamma_1$ has finite center, 
by Theorem~\ref{thm1} and Remark~\ref{Rem3.2}, 
there exists a quasi-dense subspace $X'$ of $X$ 
which splits as a product $\bar{Y}_1\times \bar{X}$, 
where $\Gamma_2\times A$ acts geometrically on $\bar{X}$ by restriction and 
$\Gamma_1$ acts geometrically on $\bar{Y}_1$ by projection.
By Theorem~\ref{thmBR} (or Theorem~\ref{thm1}), 
$\bar{X}$ contains a quasi-dense subspace $\bar{Z}\times\bar{Y}_2$, 
where $A$ acts geometrically on $\bar{Z}$ by restriction and 
$\Gamma_2$ acts geometrically on $\bar{Y}_2$ by projection.
Here 
$$ X\supset \bar{Y}_1\times \bar{X} \supset \bar{Y}_1\times \bar{Z}\times\bar{Y}_2.$$
Also, for the original action of $\Gamma$ on $X$, 
$A$ acts geometrically on $Z$ by restriction and 
$\Gamma_i$ acts geometrically on $\bar{Y}_i$ by projection for each $i=1,2$.

Then we can identify $Z=\bar{Z}$ 
because $Z$ and $\bar{Z}$ are the convex hulls of some orbits of $A$.
By construction, 
$Z\times Y_1\times Y_2$ and $\bar{Z}\times\bar{Y}_1 \times\bar{Y}_2$ are 
the unions of all subspaces which are parallel to $Z=\bar{Z}$.
Hence we obtain that $Y_1\times Y_2=\bar{Y}_1\times\bar{Y}_2$.
Here 
$\Gamma_i$ acts geometrically on $Y_i$ and $\bar{Y}_i$ 
by projection for each $i=1,2$.
Hence $Y_i$ and $\bar{Y}_i$ are the convex hulls $C(\Gamma_i*x_0)$ 
of some orbits of $\Gamma_i$ by the projection action, 
and we can identify $Y_i=\bar{Y}_i$ for each $i=1,2$.
\end{Remark}

\begin{Remark}\label{rem2}
Suppose that 
a group $A_1\times A_2\times\Gamma$ 
acts geometrically on a CAT(0) space $X$, 
where 
$A_i$ is isomorphic to $\Z^{n_i}$ for some $n_i$ ($i=1,2$) and 
$\Gamma$ has finite center.

By Theorem~\ref{thmBR}, 
$\Min (A_1\times A_2)$ is a quasi-dense subspace of $X$ and 
splits as $Z\times Y$, 
where $A_1\times A_2$ acts geometrically on $Z\cong \R^{n_1+n_2}$ by restriction 
and $\Gamma$ acts geometrically on $Y$ by projection.
Here
$$ X\supset \Min (A_1\times A_2) =Z\times Y.$$

On the other hand, 
by Theorem~\ref{thmBR}, 
$\Min (A_1)$ is a quasi-dense subspace of $X$ and 
splits as $\bar{Z}_1\times \bar{X}$, 
where $A_1$ acts geometrically on $\bar{Z}_1\cong \R^{n_1}$ by restriction 
and $A_2\times\Gamma$ acts geometrically on $\bar{X}$ by projection.
Also by Theorem~\ref{thmBR}, 
$\bar{X}$ contains a quasi-dense subspace 
$\Min_{\bar{X}}(A_2)=\bar{Z}_2\times \bar{Y}$, 
where there exist geometric actions of $A_2$ and $\Gamma$ on 
$\bar{Z}_2\cong \R^{n_2}$ and $\bar{Y}$, respectively.
Here 
\begin{align*} 
X &\supset \Min (A_1)=\bar{Z}_1\times \bar{X} \\
&\supset \bar{Z}_1\times\Min_{\bar{X}}(A_2)
=\bar{Z}_1\times\bar{Z}_2\times \bar{Y}.
\end{align*}

Then by the inductive proof 
of Theorem~\ref{thmBR} in \cite{BR} and \cite[p.245]{BH}, 
we can identify $Z=\bar{Z}_1\times\bar{Z}_2$ and $Y=\bar{Y}$.

Here we note that 
$Z=\bar{Z}_1\times\bar{Z}_2$ is the convex hull $C((A_1\times A_2)x_0)$
of some orbit of $A_1\times A_2$ and 
also we note that $Z\times Y=\bar{Z}_1\times\bar{Z}_2\times \bar{Y}$ 
is the union of all subspaces which are parallel to $Z=\bar{Z}_1\times\bar{Z}_2$ 
and $\Gamma$ acts geometrically on $Y=\bar{Y}$ by projection.
\end{Remark}

Now we prove Theorem~\ref{Thm}.

\begin{proof}[Proof of Theorem~\ref{Thm}]
Suppose that a group $\Gamma=\Gamma_1\times \Gamma_2$ acts 
geometrically on a CAT(0) space $X$.
Let $G_1=\Gamma'_1\times A_1$ and $G_2=\Gamma'_2\times A_2$ 
be essential subgroups 
of $\Gamma_1$ and $\Gamma_2$ respectively, 
where $\Gamma'_i$ has finite center and
$A_i$ is isomorphic to $\Z^{n_i}$ for some $n_i$ for each $i=1,2$.

Then 
$\Gamma'_1\times A_1\times \Gamma_2$ is a finite-index subgroup of $\Gamma$ 
and acts geometrically on $X$.
Since $\Gamma'_1$ has finite center, 
by Theorem~\ref{thm1} and Remark~\ref{Rem3.2}, 
$X$ contains a quasi-dense subspace $Y'_1\times X''_2$, 
where $\Gamma'_1$ acts geometrically on $Y'_1$ by projection and 
$A_1\times \Gamma_2$ acts geometrically on $X''_2$ by restriction.
By Theorem~\ref{thmBR}, 
$Y'_1\times X''_2$ contains a quasi-dense subspace 
$$Y'_1\times\Min_{X''_2}(A_1)=\Min_{Y'_1\times X''_2}(A_1)=
Y'_1\times Z_1\times X_2, $$
where $A_1$ acts geometrically on $Z_1$ by restriction and 
$\Gamma_2$ acts geometrically on $X_2$ by projection.
Here we note that $Y'_1\times Z_1\times X_2$ is 
a closed convex $G_1\times \Gamma_2$-invariant quasi-dense subspace of $X$.

Moreover $G_2=\Gamma'_2\times A_2$ is a finite-index subgroup of $\Gamma_2$.
Hence $Y'_1\times Z_1\times X_2$ contains a quasi-dense subspace 
$\Min_{Y'_1\times Z_1\times X_2}(A_2)=Y'_1\times Z_1\times Z_2\times Y'_2$, 
where $A_2$ and $\Gamma'_2$ act geometrically on $Z_2$ and $Y'_2$ respectively.

Thus
\begin{align*}
X &\supset Y'_1\times X_2''
\supset \Min_{Y'_1\times X_2''}(A_1)=Y'_1\times Z_1\times X_2 \\
&\supset \Min_{Y'_1\times Z_1\times X_2}(A_2)=Y'_1\times Z_1\times Z_2\times Y'_2.
\end{align*}
Here we note that $Y'_1\times Z_1\times Z_2\times Y'_2$ is 
a closed convex $G_1\times G_2$-invariant quasi-dense subspace of $X$.

Then we use Remark~\ref{rem2}.
Now $A_1\times \Gamma_2$ acts geometrically on $X''_2$ 
and $A_1\times A_2\times\Gamma'_2$ is 
a finite-index subgroup of $A_1\times \Gamma_2$.
Hence $Y'_1\times X_2''$ contains a quasi-dense subspace 
$$\Min_{Y'_1\times X_2''}(A_1\times A_2)=Y'_1\times\Min_{X''_2}(A_1\times A_2)
=Y'_1\times \hat{Y}'_2\times Z, $$
where $A_1\times A_2$ and $\Gamma'_2$ 
act geometrically on $Z$ and $\hat{Y}'_2$ respectively.
Then we can identify $Z=Z_1\times Z_2$ and $\hat{Y}'_2=Y'_2$ by Remark~\ref{rem2}.
Thus 
\begin{align*}
X &\supset Y'_1\times X_2''\supset \Min_{Y'_1\times X_2''}(A_1\times A_2) \\
&=Y'_1\times\Min_{X''_2}(A_1\times A_2)
=Y'_1\times Y'_2\times Z_1\times Z_2.
\end{align*}

Next we use Remark~\ref{rem1}.
Since $A_1\times A_2 \times\Gamma'_1\times\Gamma'_2$ 
acts geometrically on $X$, by Theorem~\ref{thmBR}, 
$X$ contains a quasi-dense subspace 
$\Min(A_1\times A_2)=\bar{Z}\times\bar{Y}'$, 
where $A_1\times A_2$ acts geometrically on $\bar{Z}$ by restriction and 
$\Gamma'_1\times\Gamma'_2$ acts geometrically on $\bar{Y}'$ by projection.
Since $\Gamma'_1\times\Gamma'_2$ has finite center, 
by Theorem~\ref{thm2}, 
$\bar{Z}\times\bar{Y}'$ contains a quasi-dense subspace 
$\bar{Z}\times\bar{Y}'_1\times \bar{Y}'_2$, 
where $\Gamma'_i$ acts geometrically on $\bar{Y}'_i$ for each $i=1,2$.
Then we can identify $\bar{Z}=Z$, $\bar{Y}'_1=Y'_1$ and $\bar{Y}'_2=\hat{Y}'_2$ 
by Remark~\ref{rem1}.
Thus 
\begin{align*}
X &\supset \Min(A_1\times A_2)=\bar{Z}\times\bar{Y'} \\
&\supset \bar{Z} \times \bar{Y}'_1\times \bar{Y}'_2 
= Z \times Y'_1\times \hat{Y}'_2 \\
&= Z_1 \times Z_2 \times Y'_1 \times Y'_2, 
\end{align*}
since $Z=Z_1 \times Z_2$ and $\hat{Y}'_2=Y'_2$.
Hence we obtain that 
$$ X\supset \Min(A_1\times A_2)\supset 
Z_1 \times Z_2 \times Y'_1 \times Y'_2 \eqno{\rm (i)}.$$

On the other hand, 
$\Gamma_1\times A_2\times \Gamma'_2$ 
is also a finite-index subgroup of $\Gamma$ and 
acts geometrically on $X$.
Since $\Gamma'_2$ has finite center, 
by Theorem~\ref{thm1} and Remark~\ref{Rem3.2}, 
$X$ contains a quasi-dense subspace $X''_1 \times \tilde{Y}'_2$, 
where $\Gamma_1\times A_2$ acts geometrically on $X''_1$ by restriction and 
$\Gamma'_2$ acts geometrically on $\tilde{Y}'_2$ by projection.
By Theorem~\ref{thmBR}, 
$X''_1 \times \tilde{Y}'_2$ contains a quasi-dense subspace 
$\Min_{X''_1 \times \tilde{Y}'_2}(A_2)=X_1\times \tilde{Z}_2\times \tilde{Y}'_2$, 
where $\Gamma_1$ and $A_2$ act geometrically on $X_1$ and $\tilde{Z}_2$ respectively.
Here we note that $X_1\times \tilde{Z}_2\times \tilde{Y}'_2$ is 
a closed convex $\Gamma_1\times G_2$-invariant quasi-dense subspace of $X$.

Moreover $G_1=\Gamma'_1\times A_1$ is a finite-index subgroup of $\Gamma_1$.
Hence $X_1\times \tilde{Z}_2\times \tilde{Y}'_2$ 
contains a quasi-dense subspace 
$\Min_{X_1\times \tilde{Z}_2\times \tilde{Y}'_2}(A_1)
=\tilde{Y}'_1\times \tilde{Z}_1\times \tilde{Z}_2\times \tilde{Y}'_2$.
Thus 
\begin{align*}
X &\supset X''_1 \times \tilde{Y}'_2
\supset \Min_{X''_1 \times \tilde{Y}'_2}(A_2)
=X_1\times \tilde{Z}_2\times \tilde{Y}'_2 \\
&\supset \Min_{X_1\times \tilde{Z}_2\times \tilde{Y}'_2}(A_1)
=\tilde{Y}'_1\times \tilde{Z}_1\times \tilde{Z}_2\times \tilde{Y}'_2.
\end{align*}
By the same argument as the above one, 
we obtain that 
$$ X\supset \Min(A_1\times A_2)
\supset \tilde{Z}_1\times \tilde{Z}_2 \times \tilde{Y}'_1\times \tilde{Y}'_2
\eqno{\rm (ii)}. $$

Since $Z_1 \times Z_2$ in (i) and $\tilde{Z}_1\times \tilde{Z}_2$ in (ii) 
are the convex hulls of some $A_1\times A_2$-orbits in $X$, 
we can identify $Z_1 \times Z_2=\tilde{Z}_1\times \tilde{Z}_2$.
Here $Z_i$ and $\tilde{Z}_i$ are isometric to $\R^{n_i}$ for each $i=1,2$ by construction.
Also 
$Z_1 \times Z_2 \times Y'_1 \times Y'_2$ in (i) and 
$\tilde{Z}_1\times \tilde{Z}_2 \times \tilde{Y}'_1\times \tilde{Y}'_2$ in (ii) 
are the unions of all subspaces which are parallel to 
$Z_1 \times Z_2=\tilde{Z}_1\times \tilde{Z}_2$ in $X$, 
hence, we can identify $Y'_1 \times Y'_2=\tilde{Y}'_1\times \tilde{Y}'_2$.
Moreover we obtain that $Y'_i=\tilde{Y}'_i$ for each $i=1,2$, 
because $Y'_i$ and $\tilde{Y}'_i$ are the convex hulls of the orbits of 
the action of $\Gamma'_i$ on $Y'_1 \times Y'_2=\tilde{Y}'_1\times \tilde{Y}'_2$ 
by projection.

Let $X'_1=Y'_1\times Z_1$ and $X'_2=\tilde{Y}'_2\times \tilde{Z}_2$.
Then 
$$X'_1\times X_2=Y'_1\times Z_1 \times X_2$$ 
is a closed convex $G_1\times \Gamma_2$-invariant quasi-dense subspace of $X$, and 
$$X_1\times X'_2=X_1\times \tilde{Z}_2\times \tilde{Y}'_2$$ 
is a closed convex $\Gamma_1\times G_2$-invariant quasi-dense subspace of $X$.
Also $X'_1=Y'_1\times Z_1$ is isometric to $\tilde{Y}'_1\times\tilde{Z}_1$ 
which is a quasi-dense subspace of $X_1$, and 
$X'_2=\tilde{Y}'_2\times \tilde{Z}_2$ is isometric to $Y'_2\times Z_2$ 
which is a quasi-dense subspace of $X_2$.
By construction, 
$\Gamma_i$ acts geometrically on $X_i$ by projection for each $i=1,2$.
Also we can obtain a geometric action of $G_i=\Gamma'_i\times A_i$ on $X'_i$ 
by the product of 
the geometric action of $\Gamma'_i$ on $Y'_i=\tilde{Y}'_i$ and 
the geometric action of $A_i$ on $Z_i\cong\tilde{Z}_i\cong\R^{n_i}$ 
for each $i=1,2$.

Here we obtain that 
the boundary $\partial X_i$ is homeomorphic to $\partial X'_i$ for each $i=1,2$, 
and 
$$\partial X=\partial (X'_1\times X_2)=
\partial X'_1*\partial X_2 \approx \partial X_1*\partial X_2.$$
Therefore the boundary $\partial X$ of $X$ is homeomorphic to 
the join $\partial X_1*\partial X_2$ of the boundaries of $X_1$ and $X_2$.
\end{proof}

A CAT(0) space $X$ is said to have the {\it geodesic extension property}, 
if every geodesic can be extended to a geodesic line $\R\rightarrow X$.

We obtain the following splitting theorem 
for CAT(0) spaces with the geodesic extension property 
from \cite[Lemma~II.6.20]{BH} and the argument in the proof of Theorem~\ref{Thm}.

\begin{Theorem}
Suppose that a group $\Gamma=\Gamma_1\times \Gamma_2$ acts geometrically 
on a CAT(0) space $X$ with the geodesic extension property.
Then 
$X$ splits as a product $X_1 \times X_2$ and 
there exist geometric actions of $\Gamma_1$ and $\Gamma_2$ 
on $X_1$ and $X_2$, respectively.
Moreover if $\Gamma$ has finite center, then 
the action of $\Gamma=\Gamma_1\times \Gamma_2$ on $X=X_1 \times X_2$ 
is the product action.
\end{Theorem}

\begin{proof}
Suppose that a group $\Gamma=\Gamma_1\times \Gamma_2$ acts geometrically 
on a CAT(0) space $X$ with the geodesic extension property.
Then by \cite[Lemma~II.6.20]{BH}, in the proof of Theorem~\ref{Thm}, 
we obtain that 
$X=X'_1 \times X_2$ where $X'_1=Y'_1\times Z_1$, 
and $X=X_1\times X'_2=\tilde{Y}'_1\times \tilde{Z}_1\times X_2'$ 
where $X'_2=\tilde{Y}'_2\times \tilde{Z}_2$ 
and $X_1=\tilde{Y}'_1\times \tilde{Z}_1$.
Here we note that 
$X_1=\tilde{Y}'_1\times \tilde{Z}_1$ is isometric to 
$X'_1=Y'_1\times Z_1$.
Hence there exists a geometric action of $\Gamma_1$ on $X'_1$.
Therefore 
$X$ splits as a product $X'_1 \times X_2$ and 
there exist geometric actions of $\Gamma_1$ and $\Gamma_2$ 
on $X'_1$ and $X_2$, respectively.

If $\Gamma$ has finite center, 
then we obtain that the action splits 
from Theorem~\ref{thm2} and \cite[Lemma~II.6.20]{BH}.
\end{proof}

\section{On products of CAT(0) groups}

A group $\Gamma$ is called a {\it CAT(0) group}, 
if $\Gamma$ acts geometrically on some CAT(0) space $X$.
A CAT(0) group $\Gamma$ is said to be {\it rigid}, 
if $\Gamma$ determines the boundary up to homeomorphisms 
of a CAT(0) space on which $\Gamma$ acts geometrically.
It is known that there exists a non-rigid CAT(0) group (\cite{CK}).

On products of CAT(0) groups, 
we obtain the following theorem from Theorem~\ref{Thm}.

\begin{Theorem}
$\Gamma_1$ and $\Gamma_2$ are CAT(0) groups if and only if 
$\Gamma_1\times \Gamma_2$ is a CAT(0) group.
\end{Theorem}

\begin{proof}
It is well known that 
if $\Gamma_1$ and $\Gamma_2$ are CAT(0) groups then 
$\Gamma_1\times \Gamma_2$ is also a CAT(0) group.
Indeed if $\Gamma_1$ and $\Gamma_2$ act geometrically 
on some CAT(0) spaces $X_1$ and $X_2$ respectively, 
then we can obtain the product geometric action of 
$\Gamma_1\times \Gamma_2$ on the CAT(0) space $X_1\times X_2$.

Conversely, 
if $\Gamma_1\times \Gamma_2$ is a CAT(0) group, 
then $\Gamma_1\times \Gamma_2$ acts geometrically on some CAT(0) space $X$ 
and by Theorem~\ref{Thm}~(4), we obtain that 
$\Gamma_1$ and $\Gamma_2$ are CAT(0) groups, 
because $X_1$ and $X_2$ in Theorem~\ref{Thm} are CAT(0) spaces.
\end{proof}

On rigidity of products of rigid CAT(0) groups, 
we also obtain the following theorem from Theorem~\ref{Thm}.

\begin{Theorem}\label{thm4}
If $\Gamma_1$ and $\Gamma_2$ are rigid CAT(0) groups, 
then so is $\Gamma_1\times \Gamma_2$, and 
the boundary $\partial(\Gamma_1\times \Gamma_2)$ 
is homeomorphic to the join $\partial\Gamma_1*\partial\Gamma_2$ 
of the boundaries of $\Gamma_1$ and $\Gamma_2$.
\end{Theorem}

\begin{proof}
Let $\Gamma_1$ and $\Gamma_2$ be rigid CAT(0) groups.
Suppose that $\Gamma_1\times\Gamma_2$ acts geometrically on a CAT(0) space $X$.
Then by Theorem~\ref{Thm}, 
there exist closed convex subspaces $X_1$ and $X_2$ in $X$ such that 
the boundary $\partial X$ is homeomorphic to the join $\partial X_1*\partial X_2$ 
and 
$\Gamma_1$ and $\Gamma_2$ act geometrically on $X_1$ and $X_2$ respectively.
Thus $\Gamma_1\times\Gamma_2$ is a rigid CAT(0) group 
whose boundary is homeomorphic to 
the join $\partial\Gamma_1*\partial\Gamma_2$ 
of the boundaries of $\Gamma_1$ and $\Gamma_2$.
\end{proof}

%

%

\begin{thebibliography}{[10]}
%
\bibitem {AB}
S.~Adams and W.~Ballmann, 
{\it Amenable isometry groups of Hadamard spaces}, 
Math.\ Ann.\ 312 (1998), 183--195.
%
\bibitem {BR}
P.~Bowers and K.~Ruane, 
{\it Boundaries of nonpositively curved groups of the form $G\times \Z^n$}, 
Glasgow Math.\ J.\ 38 (1996), 177--189.
%
\bibitem {BH}
M.~R.~Bridson and A.~Haefliger, 
{\it Metric spaces of non-positive curvature}, 
Springer-Verlag, Berlin, 1999.
%
\bibitem {BS}
M.~Burger and V.~Schroeder, 
{\it Amenable groups and stabilizers of measures on the boundary of a Hadamard manifold}, 
Math.\ Ann.\ 276 (1987), 505--514.
%
\bibitem {CK}
C.~B.~Croke and B.~Kleiner, 
{\it Spaces with nonpositive curvature and their ideal boundaries}, 
Topology 39 (2000), 549--556.
%
\bibitem {GH}
E.~Ghys and P.~de~la~Harpe (ed), 
{\it Sur les Groups Hyperboliques d'apres Mikhael Gromov}, 
Progr.\ Math.\ vol.\ 83, Birkh\"auser, Boston MA, 1990.
%
\bibitem {GW}
D.~Gromoll and J.~Wolf, 
{\it Some relations between the metric structure and the algebraic structure 
of the fundamental group in manifolds of non-positive curvature}, 
Bull.\ Amer.\ Math.\ Soc.\ 77 (1971), 545--552.
%
\bibitem {G}
M.~Gromov, {\it Hyperbolic groups}, 
Essays in group theory (S.~M.~Gersten, ed.), 
M.S.R.I. Publ.\ 8, 1987, pp.\ 75-264.
%
\bibitem {Ho}
T.~Hosaka, 
{\it A splitting theorem for CAT(0) spaces with 
the geodesic extension property}, 
Tsukuba J.\ Math.\ 27 (2003), 289--293.
%
\bibitem {JY}
J.~Jost and S.~T.~Yau, 
{\it Harmonic maps and rigidity theorems for spaces of nonpositive curvature}, 
Comm.\ Anal.\ Geom.\ 7 (1999), 681--694.
%
\bibitem {LY}
H.~B.~Lawson and S.~T.~Yau, 
{\it Compact manifolds of nonpositive curvature}, 
J.\ Differential Geom.\ 7 (1972), 211--228.
%
\bibitem {Mo}
N.~Monod, 
{\it Superrigidity for irreducible lattices and geometric splitting}, 
J.\ Amer.\ Math.\ Soc.\ 19 (2006), 781--814.
%
\bibitem {R}
K.~Ruane, 
{\it Boundaries of CAT(0) groups of the form $\Gamma=G\times H$}, 
Topology Appl.\ 92 (1999), 131--152.
%
\bibitem {S}
V.~Schroeder, 
{\it A splitting theorem for spaces of nonpositive curvature}, 
Invent.\ Math.\ 79 (1985), 323--327.
%
\end{thebibliography}
\end{document}